\documentclass[11pt]{article}
\usepackage{amssymb,amsthm}
\usepackage{amsmath}
\usepackage{amscd}
\usepackage{epsf,graphicx}
\usepackage[all]{xy}
\newtheorem{thm}{Theorem}
\newtheorem{cor}[thm]{Corollary}
\newtheorem{prop}[thm]{Proposition} 

\newtheorem{defn}[thm]{Definition}
\newtheorem{exmp}[thm]{Example}

\newcommand{\des}{\displaystyle}
\newcommand{\re}{{\rm{Re}}}

\begin{document}
\setlength{\baselineskip}{16pt}
\title{On hypergeometric functions and $k$-Pochhammer symbol}
\author{Rafael D\'\i az\thanks{Instituto Venezolano de Investigaciones Cient\'\i ficas
(IVIC). \texttt{radiaz@ivic.ve}}\ \ and Eddy
Pariguan\thanks{Universidad Central de Venezuela (UCV).
\texttt{eddyp@euler.ciens.ucv.ve}}} \maketitle

\begin{abstract}
We introduce the $k$-generalized gamma function $\Gamma_k$, beta
function $B_k$, and Pochhammer $k$-symbol $(x)_{n,k}$. We prove
several identities generalizing those satisfied by the classical
gamma function, beta function and Pochhammer symbol. We provided
integral representation for the  $\Gamma_k$ and $B_k$ functions.
\end{abstract}

\section{Introduction}
The main goal of this paper is to introduce  the $k$-gamma function
$\Gamma_k$ which is a one parameter deformation of the classical
gamma function such that $\Gamma_k\rightarrow \Gamma$ as
$k\rightarrow 1$. Our motivation to introduce $\Gamma_k$ comes from
the repeated appearance of expressions of the form
\begin{equation}\label{pocha}
x(x+k)(x+2k)\dots (x+(n-1)k)\end{equation} in a variety of
contexts, such as, the combinatorics of creation and annihilation
operators \cite{DP2}, \cite{DP1} and the perturbative computation
of Feynman integrals, see \cite{Del}. The function of variable $x$
given by formula (\ref{pocha})  will be denoted by $(x)_{n,k}$,
and will be called  the Pochhammer $k$-symbol. Setting $k=1$ one
obtains the usual Pochhammer symbol $(x)_n$, also known as the
raising factorial \cite{GCR},\ \cite{KW}. It is in principle
possible to study the Pochhammer $k$-symbol using the gamma
function, just as it is done for the case $k=1$, however one of
the main purposes of this paper is to show that it is most natural
to relate the Pochhammer $k$-symbol to the $k$-gamma function
$\Gamma_k$  to be introduce in section \ref{poga}. $\Gamma_k$ is
given by the formula
$$\Gamma_k(x)=\lim_{n\rightarrow
\infty}{\displaystyle \frac{n!
k^{n}(nk)^{{\frac{x}{k}}-1}}{(x)_{n,k}}},\ \ \ \ \ k>0,\ \ x\in
\mathbb{C}\smallsetminus k\mathbb{Z}^{-}.$$ The function $\Gamma_k$
restricted to $(0,\infty)$ is characterized by the following
properties $1)\ \ \Gamma_k(x+k)=x\Gamma_k(x)$,\ \ $2)\ \
\Gamma_k(k)=1$ and $3)\ \ \Gamma_k(x)\ \mbox{is logaritmically
convex}$.  Notice that the characterization above is indeed a
generalization of the Bohr-Mollerup theorem \cite{Con}. Just as for
the usual $\Gamma$ the function $\Gamma_k$ admits an infinite
product expression given by
\begin{equation}\frac{1}{\Gamma_k(x)}={\displaystyle{
xk^{-\frac{x}{\mbox{}k}}e^{\frac{x}{k}\gamma}\prod_{n=1}^{\infty}\left
( \left ( 1+\frac{x}{nk} \right)e^{-\frac{\!\!\!
x}{nk}}\right)}}.\end{equation} For  $\re(x)>0$, the function
$\Gamma_k$ is given by the integral
$${\des \Gamma_k(x)=\int_{0}^{\infty}
t^{x-1}e^{-\frac{t^{k}}{\!\!\!k}}dt.}$$ We deduce from the steepest
descent theorem a $k$-generalization of the famous Stirling's
formula
$$\Gamma_k(x+1)={\displaystyle
(2\pi)^{\frac{1}{2}}(kx)^{-\frac{1}{2}}x^{\frac{x+1}{k}}e^{-\frac{x}{k}}+O\left(\frac{1}{x}\right)},
\ \ \ \ \mbox{for} \ \ \ x\in \mathbb{R}^{+}.$$ It is  an
interesting problem to understand how the function $\Gamma_k$
changes as the parameter $k$ varies. Theorem \ref{pde} on section
\ref{poga} shows that the function $\psi(k,x)=\log \Gamma_{k}(x)$
is a solution of the non-linear partial differential equation
$$-kx^{2}\partial^{2}_{x}\psi
+ k^{3}\partial^{2}_{k}\psi + 2k^{2}\partial_{k}\psi=-x(k+1).$$ In
the last section of this article we study hypergeometric functions
from the point of view of the Pochhammer $k$-symbol. We
$k$-generalize some well-known identities for hypergeometric
functions such as: for any $a\in\mathbb{C}^{p}$,
$k\in(\mathbb{R}^{+})^{p}$, $s\in(\mathbb{R}^{+})^{q}$,
$b=(b_1,\dots,b_q)\in\mathbb{C}^{q}$ such that
$b_i\in\mathbb{C}\smallsetminus s_i \mathbb{Z}^{-}$ the following
identity holds
\begin{equation}
F(a,k,b,s)(x)={\des \prod_{j=1}^{p+1}\frac{1}{\Gamma_{k_j}(a_j)}
\int_{(\mathbb{R}^{+})^{p+1}}\prod_{j=1}^{p+1}
e^{\!\!-\frac{t_j^{k_j}}{k_j}} t_j^{a_j-1}
\left(\sum_{n=0}^{\infty} \frac{\!\!\!\!1}{(b)_{n,s}} \frac{(x
t_1^{k_1}\dots t_{p+1}^{k_{p+1}})^{n}}{n!}\right)
 dt},
\end{equation}
where $(b)_{n,s}=(b_1)_{n,s_1}\dots (b_q)_{n,s_q}$, $dt=dt_1\dots
dt_{p+1}$, $p\leq q$, $\re(a_j)> 0$ for all $1\leq j\leq p+1$, and
term-by-term integration is permitted. Our final result  Theorem
\ref{comint} provides combinatorial interpretation in terms of
planar forest for the coefficients of hypergeometric functions.

\section{Pochhammer $k$-symbol and $k$-gamma function}\label{poga}

In this section we present the definition of the Pochhammer
$k$-symbol and introduce the $k$-analogue of  the gamma function.
We provided representations for the $\Gamma_k$ function in term of
limits, integrals, recursive formulae, and infinite products, as
well as a generalization of the Stirling's formula.
\begin{defn}
Let  $x\in \mathbb{C}$, $k\in \mathbb{R}$ and
$n\in\mathbb{N}^{+}$, the Pochhammer $k$-symbol is given by
$$(x)_{n,k}=x(x+k)(x+2k)\dots (x+(n-1)k).$$
\end{defn}
Given  $s,n\in \mathbb{N}$ with $0\leq s \leq n$, the $s$-th
elementary symmetric function  ${\des \sum_{1\leq
i_1<\dots<i_s\leq n} x_{i_1}\dots x_{i_s}}$ on variables
$x_1,\dots,x_n$ is denoted by $e_{s}^{n}(x_1,\dots,x_n)$. Part
$(1)$ of the next proposition provides a formula for the
Pochhammer $k$-symbol in terms of the elementary symmetric
functions.

\begin{prop}
The following identities hold
\begin{enumerate}
\item{${\des (x)_{n,k}= \sum_{s=0}^{n-1}
e_{s}^{n-1}(1,2,\dots,n-1)k^{s} x^{n-s}}$.}
\item{${\des{\frac{\partial}{\partial
k}}(x)_{n,k}=\sum_{s=1}^{n-1} s(x)_{s,k} (x+(s+1)k)_{n-1-s,k}}$.}

\end{enumerate}
\end{prop}

\begin{proof} Part $(1)$ follows by induction on $n$, using the well-known
identity for elementary symmetric functions
$$e_{s}^{n-1}(x_1,\dots,x_{n-1})+n e_{s-1}^{n-1}(x_1,\dots,x_{n-1})=e_{s}^{n}(x_1,\dots,x_n).$$
Part $(2)$ follows using the logarithmic derivative.
\end{proof}

\begin{defn}\label{gengamm}
For $k>0$, the $k$-gamma function $\Gamma_k$ is given by
$$\Gamma_k(x)=\lim_{n\rightarrow \infty}{\displaystyle \frac{n!
k^{n}(nk)^{{\frac{x}{k}}-1}}{(x)_{n,k}}},\ \ \ x\in
\mathbb{C}\smallsetminus k\mathbb{Z}^{-}.$$
\end{defn}


\begin{prop}\label{impor}
Given $x\in\mathbb{C}\smallsetminus k\mathbb{Z}^{-}$, $k,s>0$ and
$n\in\mathbb{N}^{+}$, the following identity holds
\begin{enumerate}
\item{$\displaystyle{ (x)_{n,s}=\left(\frac{s}{k}\right)^{n}
\left( \frac{kx}{s}\right)_{n,k}}$.}
\item{$\displaystyle{\Gamma_{s}(x)=\left(
\frac{s}{k}\right)^{{\frac{x}{s}}-1} \Gamma_k
\left(\frac{kx}{s}\right).}$}
\end{enumerate}
\end{prop}

\begin{prop}
For $x\in \mathbb{C}$, $\re(x)>0$, we have ${\des
\Gamma_k(x)=\int_{0}^{\infty}
t^{x-1}e^{-\frac{t^{k}}{\!\!\!k}}dt.}$
\end{prop}
\begin{proof}
By Definition {\ref{gengamm}}
$$\Gamma_k(x)=\int_{0}^{\infty} t^{x-1}e^{-\frac{t^{k}}{\!\!\! k}}dt=
{\displaystyle \lim_{n\rightarrow
\infty}\int_{0}^{(nk)^{\frac{1}{k}}}\left(1-\frac{t^{k}}{nk}\right)^{n}t^{x-1}dt}.$$
Let $A_{n,i}(x), \ \ i=0,\dots, n$,  be given by
$A_{n,i}(x)={\displaystyle\int_{0}^{(nk)^{\frac{1}{k}}}\left(1-\frac{t^{k}}{nk}\right)^{i}t^{x-1}dt}.$\\

The following recursive formula is proven using integration by
parts
$$A_{n,i}(x)=\frac{i}{nx}A_{n,i-1}(x+k).$$
Also,
$$A_{n,0}(x)= {\displaystyle \int_{0}^{(nk)^{\frac{1}{k}}}t^{x-1}dt=\frac{(nk)^{\frac{x}{k}}}{x}}.$$
Therefore,
$$A_{n,n}(x)={\displaystyle \frac{n!
k^{n}(nk)^{{\frac{x}{k}}-1}}{(x)_{n,k}\left(1+\frac{x}{nk}\right)}},$$
and
$$\Gamma_k(x)=\lim_{n\rightarrow \infty}A_{n,n}(x)=\lim_{n\rightarrow \infty}{\displaystyle \frac{n!
k^{n}(nk)^{{\frac{x}{k}}-1}}{(x)_{n,k}}}. $$
\end{proof}
Notice that the case $k=2$ is of particular interest since ${\des
\Gamma_2(x)=\int_{0}^{\infty} t^{x-1}e^{-\frac{t^{2}}{\!\!\!
2}}dt}$ is the Gaussian integral.
\begin{prop}
The $k$-gamma function $\Gamma_k(x)$ satisfies the following
properties
\begin{enumerate}
\item{$\Gamma_k(x+k)=x\Gamma_k(x).$} \item{${\displaystyle
(x)_{n,k}=\frac{\Gamma_k(x+nk)}{\Gamma_k(x)}}$.}
\item{$\Gamma_{k}(k)=1$.} \item{$\Gamma_k(x)\ \ \mbox{is
logarithmically convex, for}\ x \in\mathbb{R}.$}
\item{${\displaystyle\Gamma_{k}(x)= a^{\frac{x}{k}}\int_0^{\infty}
t^{x-1} e^{-\frac{t^{k}}{\!\!\! k}a}dt}, \ \ \mbox{for}\ \
a\in\mathbb{R}$.}
 \item{${\displaystyle \frac{1}{\Gamma_k(x)}}={\displaystyle
xk^{-\frac{x}{\mbox{}k}}e^{\frac{x}{k}\gamma}\prod_{n=1}^{\infty}\left
( \left ( 1+\frac{x}{nk} \right)e^{-\frac{\!\!\! x}{nk}}\right)}$
where $\gamma={\displaystyle\lim_{n \rightarrow
\infty}(1+\cdots+\frac{1}{n}-\log(n)).}$} \item{${\displaystyle
\Gamma_k(x)\Gamma_k(k-x)= \frac{\pi}{\sin\left(\frac{\pi
x}{k}\right)}}.$}
\end{enumerate}
\end{prop}

\begin{proof}
Properties $2)$, $3)$ and $5)$ follow directly from definition.
Property $4)$ is Corollary \ref{lmc} below. $1)$, $6)$ and $7)$
follows from ${\displaystyle
\Gamma_k(x)=k^{\frac{x}{k}-1}\Gamma\left( \frac{x}{k}\right)}$.

\end{proof} Our next result is a generalization of the
Bohr-Mollerup theorem.
\begin{thm}
Let $f(x)$ be a positive valued function defined on $(0,\infty)$.
Assume that $f(k)=1$, $f(x+k)=xf(x)$ and $f$ is logarithmically
convex, then $f(x)=\Gamma_k(x)$, for all $x\in(0,\infty)$.
\end{thm}
\begin{proof}
 Identity $f(x)=\Gamma_k(x)$ holds if and only if ${\displaystyle
\lim_{n\rightarrow \infty} \frac{(x)_{n,k}f(x)}{n!
k^{n}(nk)^{{\frac{x}{k}}-1}}}=1$. Equivalently,
$$\lim_{n\rightarrow \infty}\log \left (
{\displaystyle \frac{(x)_{n,k}}{n! k^{n}(nk)^{{\frac{x}{k}}-1}}}
\right )+\log(f(x))=0.$$ Since $f$ is logarithmically convex the
following inequality holds
$${\displaystyle\frac{1}{k}\log\left (
\frac{f(nk+k)}{f(nk)}\right ) \leq \frac{1}{x}\log \left
(\frac{f(nk+k+x)}{f(nk+k)}\right ) \leq \frac{1}{k}\log \left
(\frac{f(nk+2k)}{f(nk+k)}\right )}.$$ As $f(x+k)=xf(x)$, we have
$${\displaystyle\frac{x}{k}\log(nk)\leq \log \left (\frac{(x+nk)(x+(n-1)k)\dots xf(x)}{n!k^{n}}\right )
\leq \frac{x}{k}\log((n+1)k)}$$
$${\displaystyle\log(nk)^{\frac{x}{k}}\leq \log \left (\frac{(x+nk)(x+(n-1)k)\dots xf(x)}{n!k^{n}}\right )
\leq \log((n+1)k)^{\frac{x}{k}}}$$
$${\displaystyle 0\leq \log \left (\frac{(x+nk)(x+(n-1)k)\dots xf(x)}{(nk)^{\frac{x}{k}}n!k^{n}}\right )
\leq \log \left (\frac{(n+1)k}{nk}\right )^{\frac{x}{k}}}$$
$${\displaystyle 0\leq \lim_{n \rightarrow \infty}\log \left
(\frac{(x+nk)(x+(n-1)k)\dots
xf(x)}{(nk)^{\frac{x}{k}}n!k^{n}}\right ) \leq \lim_{n \rightarrow
\infty}\log \left (\frac{(n+1)k}{nk}\right )^{\frac{x}{k}}}.$$
Since $$\lim_{n \rightarrow \infty} \log \left
(\frac{(n+1)k}{nk}\right )^{\frac{x}{k}}= \frac{x}{k}\log(1)=0,$$
we get
$${\displaystyle 0\leq \lim_{n \rightarrow \infty}\log
\left (\frac{(x+nk)(x+(n-1)k)\dots
x}{(nk)^{\frac{x}{k}}n!k^{n}}\right ) +\log(f(x))\leq 0}.$$
Therefore, $f(x)=\Gamma_k(x)$.
\end{proof}

A proof of Theorem \ref{stee} below may be found in \cite{Etg}.

\begin{thm}\label{stee}
Assume that $f:(a,b)\longrightarrow \mathbb{R}$, with  $a,b\in
[0,\infty)$ attains a global minimum at a unique point $c\in
(a,b)$, such that $f''(c)>0$. Then one has
$${\displaystyle \int_a^{b} g(x)e^{-\frac{f(x)}{\hbar}}dx=
\hbar^{\frac{1}{2}}e^{-\frac{f(c)}{\hbar}}\sqrt{2\pi}\frac{g(c)}{\sqrt{f''(c)}}+O(\hbar)}.$$

\end{thm}

As promised in the introduction, we now provide an analogue of the
Stirling's formula for $\Gamma_k$.
\begin{thm}
For $\re(x)>0$,  the following identity holds
\begin{equation}\label{stir}
\Gamma_k(x+1)={\displaystyle
(2\pi)^{\frac{1}{2}}(kx)^{-\frac{1}{2}}x^{\frac{x+1}{k}}e^{-\frac{x}{k}}
+O\left(\frac{1}{x}\right)}.
\end{equation}
\end{thm}
\begin{proof}
Recall that ${\displaystyle \Gamma_k(x+1)=\int_{0}^{\infty}
t^{x}e^{-\frac{t^{k}}{k}}dt}.$ Consider the following change of
variables $t=x^{\frac{1}{k}}v$, we get
$${\displaystyle \frac{\Gamma_k(x+1)}{x^{\frac{x+1}{k}}}=
\int_{0}^{\infty} v^{x}e^{-\frac{{(xv)}^{k}}{\!\!\!\! k}}dv=
\int_{0}^{\infty}e^{-x(\frac{v^{k}}{k}-\log v)}dv}.$$ Let
$f(s)=\frac{s^{k}}{\!\! k}-\log (s)$.  Clearly $f'(s)=0$ if and
only if $s=1$. Also $f''(1)=k$. Using Theorem \ref{stee}, we have
$$\int_{0}^{\infty}v^{x}e^{-\frac{{(xv)}^{k}}{\!\!\!\!
k}}dv=\frac{(2\pi)^{\frac{1}{2}}}{(kx)^{\frac{1}{2}}}e^{-\frac{x}{k}}
+O\left(\frac{1}{x}\right),$$ thus
$$\Gamma_k(x+1)={\displaystyle \frac{(2\pi)^{\frac{1}{2}}}{(kx)^{\frac{1}{2}}}
x^{\frac{x+1}{k}}e^{-\frac{x}{k}}+O\left(\frac{1}{x}\right)}.$$
\end{proof}
Proposition \ref{inf} and Theorem \ref{pde} bellow provide
information on the dependence of $\Gamma_k$ on the parameter $k$.

\begin{prop}\label{inf} For $\re(x)>0$, the following identity holds
$${\displaystyle\partial_{k}\Gamma_{k}(x+1)=
\frac{1}{k^{2}}\Gamma_{k}(x+k+1)-\frac{1}{k}\int_{0}^{\infty}t^{x+k}\log(t)e^{-\frac{t^{k}}{\!\!\!k}}
dt}.$$
\end{prop}
\begin{proof}
 Follows from formula
$$\Gamma_k(x+1)=\int_{0}^{\infty} t^{x}e^{-\frac{t^{k}}{\!\!\!k}}dt.
$$
\end{proof}

\begin{thm}\label{pde} For $x>0$,  the function
$\psi(k,x)=\log \Gamma_{k}(x)$ is a solution of the non-linear
partial differential equation
$$-kx^{2}\partial^{2}_{x}\psi
+ k^{3}\partial^{2}_{k}\psi + 2k^{2}\partial_{k}\psi=-x(k+1).$$
\end{thm}
\begin{proof}
Starting from
$$\frac{1}{\Gamma_k(x)}={\displaystyle{
xk^{\frac{-x}{\mbox{}k}}e^{\frac{x}{k}\gamma}\prod_{n=1}^{\infty}\left
( \left ( 1+\frac{x}{nk} \right)e^{\frac{\!\!\!
-x}{nk}}\right)}}.$$ The following equations can be proven easily.
\begin{eqnarray*}
{\displaystyle \psi(k,x)}&=&{\des
-\log(x)+\frac{x}{k}\log(k)-\frac{x}{k}\gamma
-\sum_{n=1}^{\infty}\left(\log
\left(1+\frac{x}{nk}\right)-\frac{x}{nk}\right)}.\\
{\displaystyle
\partial_{x}\psi(k,x)}&=&{\des-\frac{1}{x}+\frac{\log(k)
-\gamma}{k}-\sum_{n=1}^{\infty}\left(\frac{1}{x+nk}-\frac{1}{nk}\right)
}.\\
{\displaystyle
\partial^{2}_{x}\psi(k,x)}&=&{\des \sum_{n=0}^{\infty}\frac{1}{(x+nk)^{2}}
}.\\
{\displaystyle \partial_{k}\psi(k,x)}&=&{\des\frac{x}{k^{2}}\left(
(1-\log k+\gamma)+\sum_{n=1}^{\infty}
\left( \frac{k}{x+nk}-\frac{1}{n}\right)\right)}.\\
{\displaystyle \partial_{k}(k^{2}\partial_{k}
\psi(k,x))}&=&{\des-\frac{x}{k}+\sum_{n=1}^{\infty}
\frac{x^{2}}{(x+nk)^{2}}}.
\end{eqnarray*}

\end{proof}
The third equation above shows
\begin{cor}\label{lmc} The $k$-gamma function
$\Gamma_{k}$ is logarithmically convex on $(0,\infty)$.
\end{cor}
We remark that the $q$-analogues of the $k$-gamma and $k$-beta
functions has been introduced in \cite{DT}.

\section{$k$-beta and $k$-zeta functions}
In this section, we introduce the $k$-beta function $B_k$ and the
$k$-zeta function $\zeta_k$. We provide explicit formulae that
relate the $k$-beta $B_k$ and $k$-gamma $\Gamma_k$, in similar
fashion to the classical case.
\begin{defn}\label{betagen}
The $k$-beta function $B_k(x,y)$ is given by the formula
$$ B_{k}(x,y)=\frac{\Gamma_k(x)\Gamma_k(y)}{\Gamma_k(x+y)},\ \ \ \ \
\re(x)>0,\ \ \ \re(y)>0.$$
\end{defn}

\begin{prop}The $k$-beta function satisfies the following
identities
\begin{enumerate}
\item{${\displaystyle
B_{k}(x,y)=\int_0^{\infty}t^{x-1}(1+t^{k})^{-\frac{x+y}{k}}dt}$.}
\item{${\displaystyle
B_{k}(x,y)=\frac{1}{k}\int_0^{1}t^{\frac{x}{k}-1}(1-t)^{\frac{y}{k}-1}dt}$.}
\item{${\displaystyle B_k(x,y)=\frac{1}{k}B\left(
\frac{x}{k},\frac{y}{k}\right)}$.}
 \item{${\displaystyle
B_k(x,y)=\frac{(x+y)}{xy}\prod_{n=0}^{\infty}\frac{nk(nk+x+y)}{(nk+x)(nk+y)}}.$}
\end{enumerate}
\end{prop}

\begin{defn}
The $k$-zeta function is given by ${\displaystyle \zeta_{k}(x,s)=
\sum_{n=0}^{\infty}\frac{1}{(x+nk)^{s}}},$\ \  for $k,x>0$ and
$s>1$.
\end{defn}
\begin{thm}
The $k$-zeta function satisfies the following identities
\begin{enumerate}
\item{$\zeta_{k}(x,2)=\partial^{2}_{x}(\log \Gamma_k(x))$.}
\item{$\partial^{2}_{x}(\partial_{s}\zeta_{k})\Big{|}_{s=0}=-\partial^{2}_{x}(\log
\Gamma_{k}(x))$.} \item{${\displaystyle
\partial_{k}^{m}\zeta_{k}(x,s)=-x (s)_m\sum_{n=0}^{\infty}
\frac{n^{m}}{(x+nk)^{m+s}}}. $}
\end{enumerate}
\end{thm}
\begin{proof} Follows from equations
\begin{eqnarray*}
{\displaystyle
\partial_{s}\zeta_{k}(x,s)\Big{|}_{s=0}}&=&\des{\sum_{n=0}^{\infty}
\log(x+nk)}.\\
{\displaystyle
\partial_x(\partial_{s}\zeta_{k}(x,s))\Big{|}_{s=0}}&=&\des{
\sum_{n=0}^{\infty}\frac{1}{(x+nk)}}.\\
{\displaystyle
\partial^{2}_x(\partial_{s}\zeta_{k}(x,s))\Big{|}_{s=0}}&=&\des{-
\sum_{n=0}^{\infty}\frac{1}{(x+nk)^{2}}}.
\end{eqnarray*}
\end{proof}

\section{Hypergeometric Functions}

In this section we strongly follow the ideas and notations of
\cite{EsFu}. We study hypergeometric functions, see \cite{EsFu}
and \cite{GM} for an introduction, from the point of view of the
Pochhammer $k$-symbol.

\begin{defn}\label{hyperg}Given $a\in\mathbb{C}^{p}$, $k\in(\mathbb{R}^{+})^{p}$, $s\in(\mathbb{R}^{+})^{q}$,
$b=(b_1,\dots,b_q)\in\mathbb{C}^{q}$ such that
$b_i\in\mathbb{C}\smallsetminus s_i \mathbb{Z}^{-}$. The
hypergeometric function $F(a,k,b,s)$ is given by the formal power
series
\begin{equation}\label{hyperg1} {\des
F(a,k,b,s)(x)=\sum_{n=0}^{\infty}\frac{(a_1)_{n,k_1}(a_2)_{n,k_2}\dots(a_p)_{n,k_p}}
{(b_1)_{n,s_1}(b_2)_{n,s_2}\dots(b_q)_{n,s_q}} \frac{x^{n}}{n!}}.
\end{equation}
\end{defn}
Given $x=(x_1,\dots,x_n)\in\mathbb{R}^{n}$, we set
$\overline{x}=x_1\dots x_n$. Using the radio test one can show
that the series (\ref{hyperg1}) converges for all $x$ if $p\leq
q$. If $p>q+1$ the series diverges, and if $p=q+1$, it converges
for all $x$ such that ${\displaystyle |x|< \frac{s_1\dots
s_q}{k_1\dots k_p}}$. Also it is easy to check that the
hypergeometric function $y(x)=F(a,k,b,s)(x)$ solves the equation
$$D(s_1D+b_1-s_1)\dots (s_qD+b_q-s_q)(y)=x(k_1D+a_1)\dots(k_pD+a_p)(y),$$
where $D=x\partial_x$.

Notice that  hypergeometric function $F(a,1,b,1)$ is given by
$$F(a,1,b,1)(x)=\des{ \sum_{n=0}^{\infty} \frac{(a_1)_n\dots (a_p)_n}{(b_1)_n\dots (b_q)_n} \frac{x^{n}}{n!}},$$
and thus agrees with the classical expression for hypergeometric
functions. We now show how to transfer from the classical notation
for hypergeometric functions to our notation using the Pochhammer
$k$-symbol.

\begin{prop}\label{trans} Given $a\in\mathbb{C}^{p}$, $k\in(\mathbb{R}^{+})^{p}$, $s\in(\mathbb{R}^{+})^{q}$,
$b=(b_1,\dots,b_q)\in\mathbb{C}^{q}$ such that
$b_i\in\mathbb{C}\smallsetminus s_i \mathbb{Z}^{-}$, the following
identity holds
$${\displaystyle F(a,k,b,s)(x)=F\left(\frac{a}{k},1,\frac{b}{s},1
\right)\left( \frac{x \overline{k}}{\overline{s}}\right)},$$ where
${\des \frac{a}{k}=\left(\frac{a_1}{ k_1},\dots,\frac {a_p}{
k_p}\right)}$, ${\des \frac{b}{s}=\left(\frac{b_1}{
s_1},\dots,\frac {b_q}{ s_q}\right)}$ and $1=(1,\dots, 1)$.
\end{prop}

\begin{proof}
 \vspace{-0.2cm}
$$
F(a,k,b,s)(x)
=\sum_{n=0}^{\infty}\frac{(a)_{n,k}}{(b)_{n,s}}\frac{x^{n}}{n!}=
\sum_{n=0}^{\infty} \frac{ (\frac{a}{k})_n}{(\frac{b}{s})_n}\left(
\frac{x k_1\dots k_p} {s_1\dots s_q}\right)^{n}\frac{1}{n!}=
F\left(\frac{a}{k},1,\frac{b}{s},1 \right)\left( \frac{x
\overline{k}}{\overline{s}}\right).$$
\end{proof}

\begin{exmp}
For any  $a\in\mathbb{C}$, $k>0$ and $|x|<\frac{1}{k}$, the
following identity holds
\begin{equation}\label{bino}
{\des \sum_{n=0}^{\infty}
\frac{(a)_{n,k}}{n!}x^{n}=(1-kx)^{-\frac{a}{k}}}.
\end{equation}

\end{exmp}

We next provide an integral representation for the hypergeometric
function $F(a,k,b,s)$.  Let us first  prove a proposition that we
will be needed to obtain the integral representation. Given
$x=(x_1,\dots,x_n)\in\mathbb{C}^{n}$ we denote $x_{\leq
i}=(x_1,\dots,x_i)$.
\begin{prop} Let $a,k,b,s$ be as in Definition \ref{hyperg}.
The following identity holds
\begin{equation}\label{irep}
{\des F(a,k,b,s)(x)=\frac{1}{\Gamma_{k_{p+1}}(a_{p+1})}
\int_{0}^{\infty} e^{-\frac{t^{k_{p+1}}}{\!\!\! k_{p+1}}}
t^{a_{p+1}-1}F(a_{\leq p},k_{\leq p},b,s)(xt^{k_{p+1}})dt }
\end{equation}
when $p\leq q$, $\re(a_{p+1})> 0$,  and term-by-term integration
is permitted.
\end{prop}

\begin{proof}

${\des \int_{0}^{\infty} e^{-\frac{t^{k_{p+1}}}{\!\!\! k_{p+1}}}
t^{a_{p+1}-1}F(a_{\leq p},k_{\leq p},b,s)(xt^{k_{p+1}})dt }=$

${\des F(a_{\leq p},k_{\leq p},b,s)(x)\int_{0}^{\infty}
e^{\!\!-\frac{t^{k_{p+1}}}{\!\!\! k_{p+1}}} t^{a_{p+1}+n
k_{p+1}-1} dt} ={\des\Gamma_{k_{p+1}}(a_{p+1})F(a,k,b,s)(x) }$

\end{proof}

\begin{thm}
For any $a,k,b,s$ be as in Definition \ref{hyperg}. The following
formula holds
\begin{equation}
F(a,k,b,s)(x)={\des \prod_{j=1}^{p+1}\frac{1}{\Gamma_{k_j}(a_j)}
\int_{(\mathbb{R}^{+})^{p+1}}\prod_{j=1}^{p+1}
e^{\!\!-\frac{t_j^{k_j}}{k_j}} t_j^{a_j-1}
\left(\sum_{n=0}^{\infty} \frac{\!\!\!\!1}{(b)_{n,s}} \frac{(x
t_1^{k_1}\dots t_{p+1}^{k_{p+1}})^{n}}{n!}\right)
 dt},
\end{equation}
where $(b)_{n,s}=(b_1)_{n,s_1}\dots (b_q)_{n,s_q}$, $dt=dt_1\dots
dt_{p+1}$, $p\leq q$, $\re(a_j)> 0$ for all $1\leq j\leq p+1$, and
term-by-term integration is permitted.
\end{thm}

\begin{proof}
Use equation (\ref{irep}) and induction on $p$.
\end{proof}

\begin{exmp}
For $k=(2,2,\dots ,2)$, the hypergeometric function ${\des
F(a,2,b,s)(x)}$ is given by
$$F(a,2,b,s)={\des \prod_{j=1}^{p+1}\frac{1}{\Gamma_{2}(a_j)}
\int_{(\mathbb{R}^{+})^{p+1}}\prod_{j=1}^{p+1}
e^{\!\!-\frac{t_j^{2}}{2}} t_j^{a_j-1} \left(\sum_{n=0}^{\infty}
\frac{\!\!\!\!1}{(b)_{n,s}} \frac{(x t_1^{2}\dots
t_{p+1}^{2})^{n}}{n!}\right)
 dt},$$
where $dt=dt_1\dots dt_n$, $(b)_{n,s}=(b_1)_{n,s_1}\dots
(b_q)_{n,s_q}$, $\re(a_j)> 0$ for all $1\leq j\leq p+1$ and
term-by-term integration is permitted

\end{exmp}

We now proceed to study the combinatorial interpretation of the
coefficient of hypergeometric functions.
\begin{defn} A planar forest  $F$ consist of the following data:
\begin{enumerate}
\item{A finite totally order set $V_r(F)=\{r_1<\dots<r_m\}$ whose
elements are called roots.} \item{A finite totally order set
$V_i(F)=\{v_1<\dots<v_n\}$ whose elements are called internal
vertices.} \item {A finite set $V_t(F)$ whose elements are called
tail vertices.} \item{A map $N:V(T)\rightarrow V(T)$.} \item{Total
order on $N^{-1}(v)$ for each $v\in V(F):=V_r(F)\sqcup
V_i(F)\sqcup V_t(F)$.}
\end{enumerate}
These data satisfies the following properties:
\begin{itemize}
\item{$N(r_j)=r_j$, for all $j=1,\dots,m$ and  $N^{k}(v)=r_j$  for
some $j=1,\dots,m$ and any $k>>1$.} \item{$N(V(F))\cap
V_t(F)=\emptyset$.} \item{For any $r_j\in V_r(F)$, there is an
unique $v\in V(F)$, $v\neq r_j$ such that $N(v)=r_j$.}
\end{itemize}
\end{defn}

\begin{figure}[ht]
\centering
\includegraphics[width=3.5in]{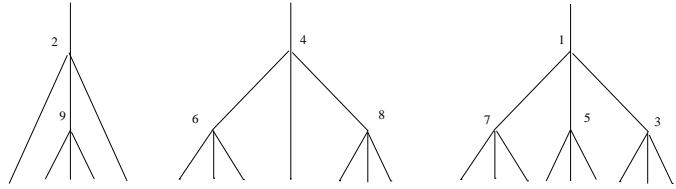}
\caption{Example of a forest in $G_{9,2}^{3}$. \label{fig:forest}}
\end{figure}

\begin{defn}a)
For any $a,k\in\mathbb{N}^{+}$, $G_{n,k}^{a}$ denotes the set of
isomorphisms classes of planar forest $F$ such that
\begin{enumerate}
\item{$V_r(F)=\{r_1<\dots<r_a\}$.}
\item{$V_i(F)=\{v_1<\dots<v_n\}$.} \item{$|N^{-1}(v_i)|=k+1$ for
all $v_i\in V_i(F)$.} \item{If $N(v_i)=v_j$, then $i<j$.}
\end{enumerate}
b) For any $a,k\in(\mathbb{N}^{+})^{p}$, we set
$G_{n,k}^{a}=G_{n,k_1}^{a_1}\times \dots \times G_{n,k_p}^{a_p}$.
\end{defn}
Figure \ref{fig:forest} provides an example of an element of
$G_{9,2}^{3}$

\begin{figure}[ht]
\centering
\includegraphics[width=2.5in]{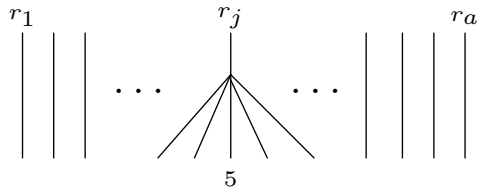}
\caption{Example of a forest in $G_{1,4}^{a}$.
\label{fig:forest2}}
\end{figure}

\begin{figure}[ht]
\centering
\includegraphics[height=1.8in]{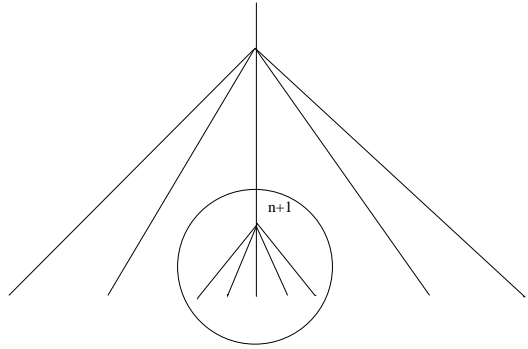}
\caption{Attaching vertex $v_{n+1}$ to a forest in $G_{n,k}^{a}$}.
\label{fig:forest1}
\end{figure}

\begin{thm}\label{comint}
Given $a,k\in(\mathbb{N}^{+})^{p}$, $b,s\in(\mathbb{N}^{+})^{q}$
and $n\in\mathbb{N}^{+}$, we have
$${\des \frac{\partial^{n}}{\partial x^{n}}F(a,k,b,s)(x)\Big{|}_{x=0}=\frac{ |G_{a,k}^{n}|}{|G_{b,s}^{n}|}.}$$
\end{thm}
\begin{proof}
 It enough to show that $(a)_{n,k}=|G_{n,k}^{a}|, \ \
\mbox{for any}\ \ a,k,n\in\mathbb{N}^{+}$. We use induction on
$n$. Since $(a)_{1,k}=a$ and $(a)_{n+1,k}=(a)_{n,k}(a+nk)$, we
have to check that $|G_{1,k}^{a}|=a$, which is obvious from Figure
\ref{fig:forest2}, and $|G_{n+1,k}^{a}|=|G_{n,k}^{a}|(a+nk)$. It
should be clear the any forest in $G_{n+1,k}^{a}$ is obtained from
a forest $F$ in $G_{n,k}^{a}$, by attaching a new vertex $v_{n+1}$
to a tail of $F$, see Figure \ref{fig:forest1}. One can prove
easily that $|V_{t}(F)|=a+nk$,  for all $F\in G_{n,k}^{a}$.
Therefore $|G_{n+1,k}^{a}|=|G_{n,k}^{a}|(a+nk)$.
\end{proof}


\newpage

\bibliographystyle{amsplain}
\bibliography{gamma}

\providecommand{\bysame}{\leavevmode\hbox to3em{\hrulefill}\thinspace}
\providecommand{\MR}{\relax\ifhmode\unskip\space\fi MR }
\providecommand{\MRhref}[2]{%
  \href{http://www.ams.org/mathscinet-getitem?mr=#1}{#2}
}
\providecommand{\href}[2]{#2}
\begin{thebibliography}{10}

\bibitem{EsFu}
Andrews, Askey, and Roy, \emph{Special functions}, Cambridge University Press,
  1999.

\bibitem{Con}
John~B. Conway, \emph{Functions of one complex variable}, {S}pringer-{V}elarg,
  {S}econd edition, New York, 1978.

\bibitem{Del}
P.~Deligne, P.~Etingof, D.~Freed, L.~Jeffrey, D.~Kazhdan, J.~Morgan,
  D.~Morrison, and E.~Witten, \emph{Quantum fields and strings: {A} course for
  mathematicians}, vol.~1, American mathematical society, 1999.

\bibitem{DP2}
Rafael D\'{\i}az and Eddy Pariguan, \emph{Quantum symmetric functions},
  math.QA/0312494, To appear in Communications in Algebra.

\bibitem{DP1}
\bysame, \emph{Symmetric quantum {Weyl} algebras}, Annales Mathematiques Blaise
  Pascal (2004), no.~11, 187--203.

\bibitem{DT}
Rafael D\'{\i}az and Carolina Teruel, \emph{q,k-generalized gamma and beta
  functions}, Journal on Nonlinear Mathematical Physics \textbf{12} (2005),
  no.~1, 118--134.

\bibitem{Etg}
Pavel Etingof, \emph{Mathematical ideas and notions of quantum field theory},
  Preprint.

\bibitem{GM}
George Gasper and Mizan Rahman, \emph{Basic hypergeometric series}, Cambridge
  University Press, New York, 1990.

\bibitem{GCR}
S.~A. Joni, G.~C Rota, and B.~Sagan, \emph{From sets to functions: {T}hree
  elementery examples}, Discrete Mathematics (1981), no.~37, 193--202.

\bibitem{KW}
K.H. Wehrhahn, \emph{Combinatorics. {A}n introduction}, vol.~2, Carslaw
  Publications, Australia, 1990.

\end{thebibliography}

\end{document}